\documentclass[12pt]{article}
\usepackage[]{amsmath,amssymb}
\usepackage{amscd}
\pdfoutput=1
\usepackage{setspace}
\usepackage{latexsym}
\usepackage{amsthm}
\usepackage{url}
\usepackage{hyperref}

\usepackage{amsfonts,  eucal} 
\usepackage[dvips]{pict2e}
\usepackage[dvipdfmx]{graphics}
\usepackage{svg}
\usepackage{float}
\usepackage{comment}
\usepackage{cite}

\newtheorem{definition}{Definition}[section]
\newtheorem{theorem}[definition]{Theorem}
\newtheorem{lemma}[definition]{Lemma}
\newtheorem{corollary}[definition]{Corollary}
\newtheorem{remark}[definition]{Remark}
\newtheorem{example}[definition]{Example}

\newcommand{\N}{\mathbb{N}}
\newcommand{\F}{\mathbb{F}}
\newcommand{\End}{\mathrm{End}}

\begin{document}
\title{\bf Concrete Billiard Arrays of Polynomial Type and Leonard Systems
}
\author{
Jimmy Vineyard
}
\date{August 28, 2024}

\maketitle
\begin{abstract} 

Let $d$ denote a nonnegative integer and let $\F$ denote a field. Let $V$ denote a $(d+1)$-dimensional vector space over $\F$. Given an ordering  $\{\theta_i\}_{i=0}^d$ of the eigenvalues of a multiplicity-free linear map $A: V \to V$, we construct a Concrete Billiard Array $\mathcal{L}$ with the property that for $0 \leq i \leq d$, the $i^{\rm th}$ vector on its bottom border is in the $\theta_i$-eigenspace of $A$. The Concrete Billiard Array $\mathcal{L}$ is said to have polynomial type. We also show the following. Assume that there exists a Leonard system $\Phi=(A;\{E_i\}_{i=0}^d;A^*;\{E_i^*\}_{i=0}^d)$ where $E_i$ is the primitive idempotent of $A$ corresponding to $\theta_i$ for $0 \leq i \leq d$. Then, we show that after a suitable normalization, the left (resp. right) boundary of $\mathcal{L}$ corresponds to the $\Phi$-split (resp. $\Phi^{\Downarrow}$-split) decomposition of $V$.

\bigskip


 \end{abstract}

\pagebreak

\section{Introduction} \label{sec:introduction}
Billiard Arrays were first introduced in \cite{Billiard_Terwilliger_2014}. In \cite{Billiard_Terwilliger_2014}, examples of Billiard Arrays were constructed using irreducible representations of  $U_q(\mathfrak{sl}_2)$ in the equitable presentation. This presentation was introduced in \cite{quantum_2006}. In \cite{Yang_2016}, it was shown that Billiard Arrays correspond to a special class of upper triangular matrices.

The concepts of Leonard pairs and Leonard systems were first introduced by Terwilliger in \cite{lp-original} in order to clarify a theorem of D. Leonard \cite{D_Leonard} concerning the orthogonal polynomials of $q$-Racah type and their relatives.

Our goal in this paper is to show how Billiard Arrays come up naturally in the context of Leonard pairs and Leonard systems. Our idea is summarized as follows. Let $V$ denote a nonzero finite dimensional vector space with dimension $d+1$.
Given an ordering  $\{\theta_i\}_{i=0}^d$ of the eigenvalues of a multiplicity-free linear map $A: V \to V$, we will construct a Concrete Billiard Array $\mathcal{L}$ with the property that for $0 \leq i \leq d$, the $i^{\rm th}$ vector on its bottom border is in the $\theta_i$-eigenspace of $A$. The Concrete Billiard Array $\mathcal{L}$ is said to have polynomial type.
We will then show the following. Suppose there exists a Leonard system
$$\Phi = (A;\{E_i\}_{i=0}^d;A^*;\{E_i^*\}_{i=0}^d),$$ 
where $E_i$ is the primitive idempotent of $A$ corresponding to $\theta_i$ for $0 \leq i \leq d$. Then, we show that after a suitable normalization, the left (resp. right) boundary of $\mathcal{L}$ corresponds to the $\Phi$-split (resp. $\Phi^{\Downarrow}$-split) decomposition of $V$.

The paper is organized as follows. Section 2 contains preliminaries. Sections 3 and 4 contain background information about Billiard Arrays. In Section 5, we introduce the Concrete Billiard Arrays of Polynomial Type. In Section 6, we use a Leonard system to construct a Concrete Billiard Array of Polynomial Type.

\section{Preliminaries} \label{sec:prelim}
Throughout this paper, we will use the following concepts and notation. Let $\F$ be a field. Let $\N = \{0,1,...\}$ be the set of natural numbers. Let $d \in \N$ and let $V$ be a $(d+1)$-dimensional $\F$-vector space. By a \textit{decomposition} of $V$, we mean a sequence $\{W_i\}_{i=0}^d$ consisting of 1-dimensional subspaces of $V$ such that the sum $V = \sum_{i=0}^d W_i $ is direct. Let $\End(V)$ denote the $\F$-algebra consisting of the $\F$-linear maps from $V$ to $V$. For $A \in \End(V)$, we say that $A$ is \textit{diagonalizable} whenever $V$ is spanned by the eigenspaces of $A$. For $A \in \End(V)$, we say that $A$ is \textit{multiplicity-free} whenever $A$ is diagonalizable and every eigenspace of $A$ has dimension 1. If $A$ is multiplicity-free then its eigenspaces form a decomposition of $V$.

\section{Billiard Arrays}\label{sec:BA}
In this section, we review the concept of a Billiard Array and a Concrete Billiard Array.
\begin{definition}\rm (\hspace{1sp}\cite[Defintion 3.1]{Billiard_Terwilliger_2014})
Let $\Delta_d$ be the subset of $\mathbb{R}^3$ consisting of all triples of natural numbers that sum to $d$, that is,
\begin{equation*}
\Delta_d = \lbrace (r,s,t) \; |\;  r,s,t \in \mathbb N, \;\;r+s+t=d\rbrace.
\end{equation*}
\noindent
An element in $\Delta_d$ is called a {\it location}. 
\end{definition}
\begin{remark} \rm \label{rmk:picture}
    We may arrange the elements of $\Delta_d$ in a triangular array. For example, if $d = 3$, $\Delta_d$ can be arranged as follows:

    \vspace{.2cm}
    \begin{center}
        030 \\
        120 \quad 021 \\
        210 \quad 111 \quad 012 \\
        300 \quad 201 \quad 102 \quad 003
    \end{center}

\end{remark}

\begin{definition}\rm
    A \textit{line} in $\Delta_d$ is a subset $L$ of $\Delta_d$ such that all $\lambda = (r,s,t) \in L$ share a common coordinate.
\end{definition}
\begin{example}
    Let $d \geq 2$. The set $\{(r,s,t) \in \Delta_d \ | \ r = 2 \}$ is a line.
\end{example}

\begin{definition}\rm
    We say two locations in $\Delta_d$ are \textit{adjacent} if they have the same value in one coordinate and differ by 1 in the other two coordinates.
    
\end{definition}

\begin{example} \rm
    Choose $d\geq1$. Let $0\leq r<d$ and $0 < s \leq d$. Then, the locations $(r,s,t)$ and $(r+1,s-1,t)$ in $\Delta_d$ are adjacent.
\end{example}

\begin{definition} \rm 
    An \textit{edge} is a pair of adjacent locations.
\end{definition}

\begin{definition}\rm (\hspace{1sp}\cite[Definition 4.29]{Billiard_Terwilliger_2014}) By a 3-clique in $\Delta_d$, we mean a set of three mutually adjacent locations in $\Delta_d$. There are two kinds of 3-cliques: $\Delta$ (black) and $\nabla$ (white).
\end{definition}

\begin{example} \rm
    Referring to the array of $\Delta_3$ in Remark \ref{rmk:picture}, the set $\{(0,3,0), (1,2,0), (0,2,1)\}$ is a black clique, and the set $\{(1,2,0), (0,2,1), (1,1,1)\}$ is a white clique.
\end{example}

\begin{definition}
\label{def:ba}
\rm
(\hspace{1sp}\cite[Definition 9.1]{Billiard_Terwilliger_2014}) Let $\mathcal{P}_1 (V)$ be the set of 1-dimensional subspaces of $V$. We define a {\it Billiard Array on $V$} to be a function
$B:\Delta_d \to \mathcal P_1(V)$, $\lambda \mapsto B_\lambda$
that satisfies the following conditions:
\begin{enumerate}
\item[\rm (i)] for each line $L$ in $\Delta_d$ 
the sum $\sum_{\lambda \in L} B_\lambda$ 
 is direct;
\item[\rm (ii)] for each black 3-clique $C$ in $\Delta_d$ the sum 
 $\sum_{\lambda \in C} B_\lambda$ 
is not direct.
\end{enumerate}
We say that $B$ is {\it over} $\mathbb F$.
We call $V$ the {\it underlying vector space}.
We call $d$ the {\it diameter} of $B$.
\end{definition}

\begin{lemma} 
\label{lem:threedim2} (\hspace{1sp}\cite[Lemma 7.3]{Billiard_Terwilliger_2014})
Let $B$ be a Billiard Array on $V$. Suppose $\lambda,\mu,\nu$ are locations in $\Delta_d$ that form a black 3-clique. Then, the subspace $B_\lambda + B_\mu + B_\nu$ is equal to each of the following:
\begin{equation}
\label{eq:directsum}
B_\lambda + B_\mu, 
\qquad
B_\mu + B_\nu, 
\qquad
B_\nu + B_\lambda.
\end{equation}
Furthermore, each sum in \eqref{eq:directsum} is direct.
\end{lemma}
\begin{proof}
Since any two of $\lambda,\mu,\nu$ are collinear, each sum in {\rm (\ref{eq:directsum})} is direct by Definition \ref{def:ba} (i). The sum $B_\lambda + B_\mu + B_\nu$ is not direct by Definition \ref{def:ba} (ii). Therefore, $B_\lambda + B_\mu$ and $B_\lambda + B_\mu + B_\nu$ are both 2-dimensional subspaces of $V$, and $B_\lambda + B_\mu$ is a subspace of $B_\lambda + B_\mu + B_\nu$. The result follows.
\end{proof}

\begin{definition}
\label{def:cba} \rm (\hspace{1sp}\cite[Definition 8.1]{Billiard_Terwilliger_2014}) We define a {\it Concrete Billiard Array on $V$}
to be a function $\mathcal B: \Delta_d \to V$,
$\lambda \mapsto \mathcal B_\lambda$ that satisfies
the following conditions:
\begin{enumerate}
\item[\rm (i)] for each line $L$ in $\Delta_d$
the vectors $\lbrace \mathcal B_\lambda \rbrace_{\lambda \in L}$
are linearly independent;
\item[\rm (ii)] for each black 3-clique $C$
 in $\Delta_d$
the vectors $\lbrace \mathcal B_\lambda \rbrace_{\lambda \in C}$
are linearly dependent.
\end{enumerate}
We say that $\mathcal B$ is {\it over} $\mathbb F$.
We call $V$ the {\it underlying vector space}.
We call $d$ the {\it diameter} of $\mathcal B$.
\end{definition}

\begin{definition}
\label{def:corr}
\rm (\hspace{1sp}\cite[Definition 8.3]{Billiard_Terwilliger_2014}) Let 
 $B$ be a Billiard Array on $V$,
 and let $\mathcal B$ be a Concrete Billiard Array on $V$.
 We say that $B$ and $ \mathcal B$ {\it correspond} whenever
 $\mathcal B_\lambda$ spans $B_\lambda$ for all 
 $\lambda \in \Delta_d$.
 \end{definition}

\begin{lemma} \label{lem:cbatoba}
    Let $\mathcal{B}$ be a Concrete Billiard Array on $V$. Then, there exists a Billiard Array $B$ on $V$ such that $B$ and $\mathcal{B}$ correspond.
\end{lemma}

\begin{proof}
    We define $B$ as follows. For $\lambda \in \Delta_{d}$, $B_\lambda = \mathrm{Span}(\mathcal{B}_\lambda)$. $B$ is a Billiard Array because Definition \ref{def:cba} satisfies the conditions of Definition \ref{def:ba}.
\end{proof}

\section{Edge Labelings and Value Functions}\label{sec:VFandEL}
In this section we review the edge labeling and value function of a Billiard Array.

Throughout this section let $B$ denote a Billiard Array with diameter $d$, and let $\mathcal{B}$ denote a Concrete Billiard Array corresponding to $B$.

Let $\lambda, \mu$ denote adjacent locations in $\Delta_d$. In \cite{Billiard_Terwilliger_2014}, a bijection $\tilde B_{\lambda,\mu}: B_\lambda \to B_\mu$ is defined. The following definitions and lemmas describe the properties of this map.
\begin{definition}
\label{def:calT}
\rm (\hspace{1sp}\cite[Definition 15.1]{Billiard_Terwilliger_2014})
Let $\lambda, \mu$ denote adjacent locations in $\Delta_d$. Recall that $\mathcal B_\lambda $ is a basis for $B_\lambda$
and 
 $\mathcal B_\mu $ is a basis for $B_\mu$.
Define a scalar
$\tilde {\mathcal B}_{\lambda,\mu} \in \mathbb F$ such that
$\tilde B_{\lambda,\mu}$ sends
	  $
	  \mathcal B_\lambda \mapsto  
      \tilde {\mathcal B}_{\lambda,\mu} \mathcal B_\mu
$.
Note that 
$\tilde {\mathcal B}_{\lambda,\mu} \not=0$.
\end{definition}

\begin{lemma} 
\label{lem:dep} (\hspace{1sp}\cite[Lemma 15.6]{Billiard_Terwilliger_2014})
Let $\lambda, \mu, \nu$ denote
locations in $\Delta_d$ that form a black 3-clique.
Then
\begin{equation*}
\label{eq:ld}
\mathcal B_\lambda
+ 
\tilde {\mathcal B}_{\lambda, \mu}\mathcal B_\mu
+
\tilde {\mathcal B}_{\lambda, \nu}\mathcal B_\nu = 0.
\end{equation*}
\end{lemma}

\begin{lemma} 
\label{lem:whiteprod} (\hspace{1sp}\cite[Lemma 15.8]{Billiard_Terwilliger_2014})
Let $\lambda, \mu, \nu$ denote
locations in $\Delta_d$ that form a black 3-clique.
Then
\begin{equation*}
          \tilde {\mathcal B}_{\lambda,\mu}
          \tilde {\mathcal B}_{\mu,\nu}
          \tilde {\mathcal B}_{\nu,\lambda}
	  = 1.
\end{equation*}
\end{lemma}

\begin{lemma}
\label{lem:CBAreciprocal} (\hspace{1sp}\cite[Lemma 15.4]{Billiard_Terwilliger_2014})
Let $\lambda, \mu$ denote adjacent
locations in $\Delta_d$.
Then the scalars
$\tilde {\mathcal B}_{\lambda,\mu}$ and
$\tilde {\mathcal B}_{\mu,\lambda}$ are reciprocal.
\end{lemma}

\begin{definition}
    \label{def:B-value} (\hspace{1sp}\cite[Definition 14.9]{Billiard_Terwilliger_2014}) Let $\lambda, \mu, \nu$ denote locations in $\Delta_d$ that form a white $3$-clique. The nonzero scalar $\tilde {\mathcal B}_{\lambda,\mu}\tilde {\mathcal B}_{\mu,\nu}\tilde {\mathcal B}_{\nu,\lambda}$ is called the {\it clockwise $\beta$-value} 
(resp. {\it counterclockwise $\beta$-value})
of the clique whenever
the sequence 
 $\lambda, \mu, \nu$ runs clockwise (resp. counterclockwise)
 around the clique.
\end{definition}

We now define the notion of an edge labeling.

\begin{definition}\rm (\hspace{1sp}\cite[Defintion 16.1]{Billiard_Terwilliger_2014})
\label{def:EL}
 By an {\it edge labeling of $\Delta_d$}
we mean a function $\beta$
that assigns to each ordered pair 
$\lambda, \mu$ of adjacent locations in $\Delta_d$
a scalar $\beta_{\lambda, \mu} \in  \mathbb F$  such that:
\begin{enumerate}
\item[\rm (i)] for adjacent locations $\lambda, \mu$ in $\Delta_d$,
\begin{equation*}
\beta_{\lambda, \mu}\beta_{\mu,\lambda} = 1;
\end{equation*}
\item[\rm (ii)] for locations $\lambda, \mu, \nu$  in $\Delta_d$
that form a black 
3-clique,
\begin{equation*}
\beta_{\lambda, \mu}
\beta_{\mu, \nu}
\beta_{\nu, \lambda}
= 1.
\end{equation*}
\end{enumerate}
\end{definition}

\begin{lemma}
\label{prop:CBAtoEL} (\hspace{1sp}\cite[Lemma 16.2]{Billiard_Terwilliger_2014})
Define a function $\tilde {\mathcal B}$ that assigns to each ordered pair $\lambda, \mu$ of adjacent locations in $\Delta_d$ the scalar $\tilde {\mathcal B}_{\lambda,\mu}$ from Definition \ref{def:calT}.
Then $\tilde {\mathcal B}$ is an edge labeling of $\Delta_d$.
\end{lemma}
\begin{proof}
The function $\tilde {\mathcal B}$
satisfies the conditions (i) and (ii)  of
Definition
\ref{def:EL} by
 Lemmas
\ref{lem:CBAreciprocal} and 
\ref{lem:whiteprod} respectively.
\end{proof}

\begin{lemma} (\hspace{1sp}\cite[Lemma 16.3]{Billiard_Terwilliger_2014})
\label{lem:ELnonzero}
For adjacent locations $\lambda, \mu$ in $\Delta_d$ we have
$\beta_{\lambda,\mu}\not=0$.
\end{lemma}
\begin{proof} By Definition \ref{def:EL} (i).
\end{proof}

\begin{definition} 
\label{def:orientedbetaval}
\rm (\hspace{1sp}\cite[Definition 18.1]{Billiard_Terwilliger_2014})
Let $\beta$ denote an edge labeling of $\Delta_d$.
Let $\lambda, \mu, \nu$ denote locations in $\Delta_d$
that form a white 3-clique. Then the scalar
\begin{equation*}
\beta_{\lambda,\mu}
\beta_{\mu,\nu}
\beta_{\nu,\lambda}
\end{equation*}
is called the {\it clockwise $\beta$-value} 
(resp. {\it counterclockwise $\beta$-value})
of the clique whenever
the sequence 
 $\lambda, \mu, \nu$ runs clockwise (resp. counterclockwise)
 around the clique.
\end{definition}

\begin{lemma} (\hspace{1sp}\cite[Lemma 18.2]{Billiard_Terwilliger_2014})
Let $\beta$ denote an edge labeling of $\Delta_d$.
For each white 3-clique in $\Delta_d$, 
its clockwise $\beta$-value and counterclockwise $\beta$-value are reciprocal.
\end{lemma}
\begin{proof}
By Definition \ref{def:EL} (i) and
Definition
\ref{def:orientedbetaval}.
\end{proof}

\begin{definition}
\label{def:betavalStraight}
\rm (\hspace{1sp}\cite[Definition 18.3]{Billiard_Terwilliger_2014})
Let $\beta$ denote an edge labeling of $\Delta_d$.
 For each white 3-clique in $\Delta_d$,
by its {\it $\beta$-value}  we mean its clockwise $\beta$-value.
\end{definition}

\begin{lemma}
\label{lem:black}
\rm (\hspace{1sp}\cite[Lemma 4.32]{Billiard_Terwilliger_2014}) Assume $d\geq 2$. 
We describe a bijection from $\Delta_{d-2}$ to the
set of white 3-cliques in $\Delta_d$. The
bijection sends each
$(r,s,t) \in \Delta_{d-2}$ to the white 3-clique consisting
of the locations 
\begin{equation*}
(r,s+1,t+1), \qquad (r+1,s,t+1), \qquad (r+1,s+1,t).
\end{equation*}
\end{lemma}

\begin{definition}
\label{def:vfunction} (\hspace{1sp}\cite[Definition 14.12]{Billiard_Terwilliger_2014})
\rm By a {\it value function} on 
$\Delta_d$ we mean a function $\psi: \Delta_d \to
{\mathbb F}\backslash \{ 0 \}$.
\end{definition}

\begin{definition} \label{def:val-func-B} (\hspace{1sp}\cite[Definition 14.13]{Billiard_Terwilliger_2014})
Assume $d \geq 2$. We define a function $\hat{B} : \Delta_{d-2} \to \mathbb{F}$ as follows. Pick $(r, s, t) \in \Delta_{d-2}$. To describe the image of $(r, s, t)$ under $\hat{B}$, consider the corresponding white 3-clique in $\Delta_d$ from Lemma \ref{lem:black}. The $B$-value of this 3-clique is the image of $(r, s, t)$ under $\hat{B}$. Observe that $\hat{B}$ is a \textit{value function} on $\Delta_{d-2}$, in the sense of Definition \ref{def:vfunction}. We call $\hat{B}$ the \textit{value function} for $B$.
\end{definition}

\begin{definition}
\label{def:el}
\rm (\hspace{1sp}\cite[Definition 16.1]{Billiard_Terwilliger_2014}) Let $d\geq 2$, and let
$\beta$ be an edge labeling of $\Delta_d$.
We define a function
$\hat \beta:\Delta_{d-2} \to \mathbb F$ as follows.
Pick $(r,s,t) \in \Delta_{d-2}$. The image of $(r,s,t)$ under $\hat \beta$ is determined by considering the associated white 3-clique in $\Delta_d$ from Lemma \ref{lem:black}. The $\beta$-value assigned to this 3-clique is the image of $(r,s,t)$ under $\hat \beta$. $\hat \beta$ serves as a value function on $\Delta_{d-2}$ by Definition \ref{def:vfunction}. We refer to $\hat \beta$ as the {\it value function} for $\beta$. 
\end{definition}

The following are some main results about Billiard Arrays.

\begin{definition} \rm
    Let $B,B'$ be Billiard Arrays on the $(d+1)$-dimensional $\F$-vector spaces $V$ and $V'$ respectively. An \textit{isomorphism of Billiard Arrays} from $B$ to $B'$ is an $\F$-linear bijection $\sigma: V \to V'$ that sends $B_\lambda \mapsto B_\lambda'$ for all $\lambda \in \Delta_d$. Two Billiard Arrays are \textit{isomorphic} if there exists an isomorphism between them.
\end{definition}

\begin{lemma}
\label{lem:diagcom} (\hspace{1sp}\cite[Lemma 18.5]{Billiard_Terwilliger_2014})
Let $\mathcal B$ denote a Concrete Billiard
Array over $\mathbb F$ that has diameter $d$. Let $B$
denote the corresponding Billiard Array. The value functions for $B$ and $\tilde{\mathcal{B}}$ as described in Definitions \ref{def:val-func-B} and \ref{def:el} are the same.
\end{lemma}

\section{Concrete Billiard Arrays of Polynomial Type} \label{sec:diag-type}
In this section, we introduce the concept of a Concrete Billiard Array of Polynomial Type.

Throughout this section, let $A \in \End(V)$ be multiplicity-free. We fix an ordering $\theta_0,...,\theta_d$ of the eigenvalues of $A$. For $0 \leq i \leq d$, let $E_i$ be the element of $\End(V)$ that acts as the identity on the $\theta_i$-eigenspace of $A$, and as 0 on the other eigenspaces of $A$. The map $E_i$ is sometimes called the primitive idempotent of $A$ corresponding to $\theta_i$. It is well known that
\begin{equation*}
        E_i = \prod_{\substack{0 \leq j \leq d \\ j \neq i}} \frac{A - \theta_j I}{\theta_i - \theta_j}.
    \end{equation*}
For $0 \leq i \leq d$, pick nonzero $\mathbf{v}_i \in E_iV$. Let $\mathbf{v} = \mathbf{v}_0+\mathbf{v}_1+\cdots + \mathbf{v}_d$. Observe that for $0 \leq i \leq d$, $E_i\mathbf{v} = \mathbf{v}_i \neq 0$. Let $\xi$ denote an indeterminate.

\begin{definition} \label{def:tau-eta} \rm For $0 \leq k \leq d$, define the polynomials
    \begin{equation*}
        \tau_k(\xi) =  \prod_{i=1}^k (\xi-\theta_{i-1}), \qquad \qquad \eta_k(\xi) = \prod_{i=1}^k (\xi-\theta_{d-i+1}).
    \end{equation*}
    We interpret $\tau_0(\xi) = \eta_0(\xi) = 1$.
\end{definition}

\begin{lemma} \label{lem:eta-zero}
    For $0 \leq i,j\leq d$ the following holds:
    \begin{align*}
        \tau_j(\theta_i) &= 0 \quad \iff \quad i < j, \\
        \eta_j(\theta_i) &= 0 \quad \iff \quad i > d-j.
    \end{align*}
\end{lemma}

\begin{proof}
    We verify the first displayed line. If $i < j$, then $\xi - \theta_i$ is a factor of $\tau_j(\xi)$. If $i \geq j$, then $\xi - \theta_i$ is not a factor of $\tau_j(\xi)$, and therefore $\theta_i$ is not a root. The first displayed line is verified. The second displayed line is similarly verified.
\end{proof}

\begin{definition} \rm \label{def:leonard-array}
    For all $\lambda = (r,s,t) \in \Delta_d$, we define a vector $\mathcal{L}_\lambda \in V$ as follows:
    \begin{equation*}
        \mathcal{L}_\lambda = \eta_r(A)\tau_t(A)\mathbf{v}.
    \end{equation*}
    We define the set 
    \begin{equation*}
       \mathcal{L} = \{\mathcal{L}_\lambda \ | \ \lambda \in \Delta_d\}.
    \end{equation*}
\end{definition}

Shortly we will show that $\mathcal{L}$ is a Concrete Billiard Array. 

\begin{lemma} \label{lem:tau-eta-A}
    For $0 \leq i,j \leq d$ the following hold:
    \begin{align*}
        E_i\tau_j(A)\mathbf{v} &= \tau_j(\theta_i)\mathbf{v}_i,\\ 
        E_i \eta_j(A)\mathbf{v} &= \eta_j(\theta_i)\mathbf{v}_i.
    \end{align*}
\end{lemma}
\begin{proof}
    Use Definition \ref{def:tau-eta} and the definition of $\mathbf{v}$ along with the fact that $E_iA = \theta_iE_i$.
\end{proof}

\begin{lemma} \label{lem:line-indep}
    For $0 \leq k \leq d$ the following sets are linearly independent:
    \begin{align*}
       R_k &= \{\eta_r(A)\tau_t(A)\mathbf{v} \ |  \ (r,s,t) \in \Delta_d, r = k\},\\
       S_k &= \{\eta_r(A)\tau_t(A)\mathbf{v} \ |\  (r,s,t) \in \Delta_d, s = k\},\\
       T_k &= \{\eta_r(A)\tau_t(A)\mathbf{v} \ | \ (r,s,t) \in \Delta_d, t = k\}.
    \end{align*}
\end{lemma}

\begin{proof}
    We first show this for $S_k$. Pick scalars $\alpha_0,...,\alpha_{d-k} \in \F$ such that
    \begin{equation} \label{equ:indep-S}
        \sum_{t=0}^{d-k} \alpha_t\eta_{d-k-t}(A)\tau_t(A)\mathbf{v} = 0.
    \end{equation}
    We will show $\alpha_t = 0$ for $0 \leq t \leq d-k$. Suppose there exists at least one $t$ ($0 \leq t \leq d-k$) such that $\alpha_t \neq 0$. Define $n = \min\{t \ | \  0 \leq t \leq d-k, \ \alpha_t \neq 0\}$. Observe that
    \begin{equation*}
        \sum_{t=0}^{d-k} \alpha_t\eta_{d-k-t}(A)\tau_t(A)\mathbf{v} = \sum_{t=n}^{d-k} \alpha_t\eta_{d-k-t}(A)\tau_t(A)\mathbf{v}.
    \end{equation*}
    By Lemmas \ref{lem:eta-zero}, \ref{lem:tau-eta-A} and the definition of $\mathbf{v}$, we have
    \begin{align*} \label{equ:indep-S-2}
        E_n \sum_{t=0}^{d-k} \alpha_t\eta_{d-k-t}(A)\tau_t(A)\mathbf{v}&= E_n \sum_{t=n}^{d-k} \alpha_t\eta_{d-k-t}(A)\tau_t(A)\mathbf{v} \\
        &= E_n\alpha_n\eta_{d-k-n}(A)\tau_n(A)\mathbf{v} \\
        &= \alpha_n\eta_{d-k-n}(\theta_n)\tau_n(\theta_n)E_n\mathbf{v}\\
        &\neq 0,
    \end{align*}
    contradicting \eqref{equ:indep-S}.
    
    We have proven the result for $S_k$. The proof for $R_k$ is similar. The proof for $T_k$ is also similar, except that we invert the order of the $E_i$'s.
\end{proof}

\begin{lemma} \label{lem:black-clique-relation}
    Let $\lambda, \mu, \nu \in \Delta_d$ form a black clique, arranged as 
    $$\lambda = (r,s,t), \qquad \mu = (r+1,s-1,t), \qquad \nu = (r,s-1,t+1).$$
    Then the following equations hold.
    \begin{align*}
        \mathcal{L}_\mu = (A-\theta_{d-r}I)\mathcal{L}_\lambda, \qquad \quad
        \mathcal{L}_\nu = (A-\theta_tI)\mathcal{L}_\lambda, \qquad \quad
        \mathcal{L}_\lambda = \frac{\mathcal{L}_\mu - \mathcal{L}_\nu}{\theta_t - \theta_{d-r}}.
    \end{align*}
\end{lemma}

\begin{proof}  
    The result for $\mathcal{L}_\mu$ and $\mathcal{L}_\nu$ is immediate from Definitions \ref{def:tau-eta} and \ref{def:leonard-array}. The result for $\mathcal{L}_\lambda$ follows from the construction.
\end{proof}

We emphasize one aspect of Lemma \ref{lem:black-clique-relation}.
\begin{lemma} \label{lem:black-indep}
        The vectors $\mathcal{L}_\lambda, \mathcal{L}_\mu, \mathcal{L}_\nu$ from Lemma \ref{lem:black-clique-relation} are linearly dependent.
\end{lemma}

\begin{proof}
    Use Lemma \ref{lem:black-clique-relation}.
\end{proof}

\begin{theorem} \label{thm:lp-to-cba}
    The set $\mathcal{L}$ from Definition \ref{def:leonard-array} is a Concrete Billiard Array.
\end{theorem}
\begin{proof}
    We verify the two conditions from Definition \ref{def:cba}. Condition (i) is satisfied by Lemma \ref{lem:line-indep} and condition (ii) is satisfied by Lemma \ref{lem:black-indep}.
\end{proof}

\begin{definition} \rm \label{def:diagonal-array}
    We call the set $\mathcal{L}$ from Definition \ref{def:leonard-array} the \textit{Concrete Billiard Array of Polynomial Type} corresponding to $A$ and $\mathbf{v}$.
\end{definition}

Recall the picture from Remark \ref{rmk:picture}. In the next result, we consider a location on the bottom border.

\begin{lemma} \label{lem:eigenspace-decomp-array}
For $0 \leq i \leq d$, consider the location $\lambda_i = (d-i,0,i) $ in $ \Delta_d$. Then $\mathcal{L}_{\lambda_i} \in E_iV$.
\end{lemma}
\begin{proof}
    Use Definitions \ref{def:tau-eta} and \ref{def:leonard-array}.
\end{proof}

Our next goal is to find the value function for the Concrete Billiard Array $\mathcal{L}$.

 Recall the bijections $\tilde B_{\lambda,\mu}: B_\lambda \to B_\mu$ from Definition \ref{def:calT}. In that definition, we encountered the scalars $\tilde{\mathcal{B}}_{\lambda,\mu}$. For notational convenience, in our discussion of $\mathcal{L}$ we will refer to these scalars as $\tilde{\mathcal{L}}_{\lambda,\mu}$.
\begin{lemma} \label{lem:lp-edgelabel}
    
    Let $\lambda, \mu, \nu \in \Delta_d$ form a black clique, arranged as 
    $$\lambda = (r,s,t), \qquad \mu = (r+1,s-1,t), \qquad \nu = (r,s-1,t+1).$$
    Going counterclockwise around this clique, the scalars from Definition \ref{def:calT} are
    $$ \tilde{\mathcal{L}}_{\lambda\mu} = \frac{1}{\theta_t - \theta_{d-r}}, \qquad \tilde{\mathcal{L}}_{\mu\nu} = -1, \qquad \tilde{\mathcal{L}}_{\nu\lambda} = \theta_{d-r} - \theta_t. $$
\end{lemma}

\begin{proof}
    By Lemma \ref{lem:black-clique-relation}, we have $\mathcal{L}_\lambda = \frac{1}{\theta_t - \theta_{d-r}}(\mathcal{L}_\mu - \mathcal{L}_\nu)$. We may rearrange this to get 
    \begin{equation} \label{equ:L-sum-1}
    \mathcal{L}_\lambda + \frac{1}{\theta_{d-r}-\theta_t}\mathcal{L}_\mu + \frac{1}{\theta_t-\theta_{d-r}}\mathcal{L}_\nu = 0.
    \end{equation}
    By Lemma \ref{lem:dep}, 
    \begin{equation} \label{equ:L-sum-2}
    \mathcal L_\lambda + \tilde {\mathcal L}_{\lambda, \mu}\mathcal L_\mu + \tilde {\mathcal L}_{\lambda, \nu}\mathcal L_\nu = 0.
    \end{equation}
    Comparing \eqref{equ:L-sum-1} and \eqref{equ:L-sum-2}, by Defintion \ref{def:cba} we obtain 
    $$\tilde {\mathcal{L}}_{\lambda,\mu} = \frac{1}{\theta_{d-r}-\theta_t}, \qquad \qquad \tilde {\mathcal{L}}_{\lambda,\nu} = \frac{1}{\theta_t-\theta_{d-r}}.$$ 
    By Lemma \ref{lem:CBAreciprocal}, 
    $$\tilde {\mathcal{L}}_{\nu,\lambda} = \frac{1}{\tilde {\mathcal{L}}_{\lambda,\nu}}= \theta_t-\theta_{d-r}.$$
    By Lemma \ref{lem:whiteprod}, 
    $$\tilde {\mathcal{L}}_{\lambda,\mu}\tilde {\mathcal{L}}_{\mu,\nu}\tilde {\mathcal{L}}_{\nu,\lambda} = 1.$$
    By these comments $\tilde{\mathcal{L}}_{\mu,\nu} = -1$. 
\end{proof}

\begin{corollary} \label{cor:edge-label}
    Let $\lambda, \mu, \nu \in \Delta_d$ form a black clique. Consider the map that sends the adjacent locations $\lambda,\mu$ (resp. $\lambda, \nu$) (resp. $\mu, \nu)$ in the black clique to the scalar $\tilde{\mathcal{L}}_{\lambda,\mu}$ (resp. $\tilde{\mathcal{L}}_{\lambda,\nu}$) (resp. $\tilde{\mathcal{L}}_{\mu,\nu}$) from Lemma \ref{lem:lp-edgelabel}. This map is an edge labeling on $\mathcal{L}$.
\end{corollary}
\begin{proof}
    Use Lemma \ref{prop:CBAtoEL}.
\end{proof}

\begin{theorem} \label{thm:la-value-function}
    The clockwise value function $\Delta_{d-2} \to \F$ for $\mathcal{L}$ sends 
    \begin{equation} \label{equ:white-value-la}
        (r,s,t) \mapsto \frac{\theta_{d-r-1} - \theta_t}{\theta_{d-r} - \theta_{t+1}}.
    \end{equation}
\end{theorem}
\begin{proof}
    Let $(r,s,t) \in \Delta_{d-2}$. Consider the following locations in $\Delta_{d}$:
    \begin{equation*}
        \lambda = (r+1,s+1,t), \qquad \mu = (r, s+1, t+1), \qquad \nu = (r+1,s,t+1).
    \end{equation*}
    These locations form a white clique in $\Delta_{d}$ with the following arrangement:
    \begin{center}
       \begin{tabular}{ccc}
        $\lambda$ & & $\mu$ \\
         & $\nu$ &
    \end{tabular} 
    \end{center}
    This clique corresponds to the location $(r,s,t) \in \Delta_{d-2}$.
    By Definitions $\ref{def:orientedbetaval}$ and $\ref{def:el}$, the value function sends $(r,s,t)$ to $\tilde {\mathcal{L}}_{\lambda,\mu}\tilde {\mathcal{L}}_{\mu,\nu}\tilde {\mathcal{L}}_{\nu,\lambda}$. By Lemma \ref{lem:lp-edgelabel},
    \begin{equation*}
        \tilde {\mathcal{L}}_{\lambda,\mu} = -1, \qquad \tilde {\mathcal{L}}_{\mu,\nu} = \frac{1}{\theta_{t+1}-\theta_{d-r}}, \qquad \tilde {\mathcal{L}}_{\nu,\lambda} = \theta_{d-r-1} -\theta_t.
    \end{equation*}
    Observe that 
    $$\tilde {\mathcal{L}}_{\lambda,\mu}\tilde {\mathcal{L}}_{\mu,\nu}\tilde {\mathcal{L}}_{\nu,\lambda} = \frac{\theta_{d-r-1} - \theta_t}{\theta_{d-r} - \theta_{t+1}}.$$
    The result follows.
\end{proof}

\section{Application to Leonard Systems} \label{sec:leonard}
In this section, we use a Leonard system to obtain a Concrete Billiard Array of Polynomial Type. We begin our discussion with a review of Leonard pairs.
\begin{definition}\rm (\hspace{1sp}\cite{lp-original})
\label{def:lp-def}
By a {\it Leonard pair} on $V$,
we mean an ordered pair $A, A^* $ of elements in $ \End(V)$
that satisfies the following conditions: 
\begin{enumerate}
\item there exists a basis for $V$ with respect to which
the matrix representing $A$ is diagonal and the matrix
representing $A^*$ is irreducible tridiagonal;
\item there exists a basis for $V$ with respect to which
the matrix representing $A^*$ is diagonal and the matrix
representing $A$ is irreducible tridiagonal.

\end{enumerate}
(A tridiagonal matrix is said to be irreducible
whenever all entries on the super and sub diagonals are nonzero).

\end{definition}

\begin{lemma} \label{lem:lp-eigenvalues}(\hspace{1sp}\cite{lp-original})
    Let $A, A^*$ denote a Leonard pair on $V$. Then
    $A$ and $A^*$ are multiplicity-free.
\end{lemma}

\begin{definition} \rm (\hspace{1sp}\cite{lp-original})
\label{def:defls}
By a {\it Leonard system} on $V$, we mean a 
sequence 
\begin{equation*}
\; \Phi = (A;\{E_i\}_{i=0}^d;\,A^*;\{E_i^*\}_{i=0}^d)
\label{eq:ourstartingpt}
\end{equation*}
 that satisfies the following conditions:
\begin{enumerate}
\item $A,A^*$ are both multiplicity-free elements in $\End(V)$;
\item $\{E_i\}_{i=0}^d$ is an ordering of the primitive idempotents of $A$;
\item $\{E_i^*\}_{i=0}^d$ is an ordering of the primitive idempotents of $A^*$;
\item $E_iA^*E_j = \begin{cases}
                0, &\mathrm{if} \ \ \ \ |i - j| > 1; \\
                \neq 0 &\mathrm{if} \ \ \ \ |i-j| = 1;
            \end{cases}$
\item $E_i^*AE_j^* = \begin{cases}
                0, &\mathrm{if} \ \ \ \ |i - j| > 1; \\
                \neq 0 &\mathrm{if} \ \ \ \ |i-j| = 1.
            \end{cases}$
\end{enumerate}
We refer to $d$ as the {\it diameter} of $\Phi$, and say 
$\Phi$ is {\it over} $\F$.

\end{definition}
In order to motivate our next results, we recall a few facts about Leonard systems.

Note that the following is a Leonard system on $V$.
\begin{equation} \label{equ:phi-downarrow}
    \Phi^{\Downarrow} = (A;\{E_{d-i}\}_{i=0}^d;A^*;\{E^*_i\}_{i=0}^d).
\end{equation}

 We now recall the split decompositions of $V$.
For $0 \leq i \leq d$ we define
\begin{equation}
U_i = 
(E^*_0V + E^*_1V + \cdots + E^*_iV)\cap (E_0V + E_{i+1}V + \cdots + E_{d-i}V).
\label{eq:defui}
\end{equation}
It was shown in \cite{TD00}
that each of $U_0, U_1, \ldots, U_d$ has dimension 1, and that
\begin{equation}
V = U_0 + U_1 + \cdots + U_d \qquad \qquad (\hbox{direct sum}).
\label{eq:splitdec}
\end{equation}
Moreover, the following holds for $0 \leq i \leq d$.
\begin{eqnarray}
U_0 + U_1 + \cdots + U_i &=& E^*_0V + E^*_1V + \cdots + E^*_iV,
\label{eq:split1}
\\
U_i + U_{i+1} + \cdots + U_d &=& E_0V + E_{1}V + \cdots + E_{d-i}V
\label{eq:split2}.
\end{eqnarray}
We call $\{U_i\}_{i=0}^d$ the $\Phi$\textit{-split decomposition} of $V$. Similarly, we call $\{U^{\Downarrow}_i\}_{i=0}^d$ the $\Phi^{\Downarrow}$\textit{-split decomposition} of $V$, where the subspaces $U_i^\Downarrow$ are defined analogously to $U_i$, but with respect to $\Phi^\Downarrow$ instead of $\Phi$.
\noindent
By \cite{TD00,lp-original}, the following holds for $0 \leq i \leq d$.
\begin{align}
    U_i &= (A-\theta_d I)(A - \theta_{d-1} I)\cdots (A-\theta_{d-i+1} I )E^*_0V = \eta_i(A)E_0^*V.
\label{eq:uialt}\\
    U_i^{\Downarrow} &= (A-\theta_0 I)(A - \theta_{1} I)\cdots (A-\theta_{i-1} I )E^*_0V = \tau_i(A)E_0^*V.\label{eq:uialt2}
\end{align}
\noindent
In particular,
$$U_0 = E_0^*V, \qquad \qquad U_d = E_0V.$$

We now discuss how Leonard systems are related to Concrete Billiard Arrays of Polynomial Type.

Referring to the Leonard system in Definition \ref{eq:ourstartingpt}, the map $A$ is multiplicity-free. Therefore, by Definition \ref{def:leonard-array}, we may use $A$ to create a Concrete Billiard Array $\mathcal{L}$ of Polynomial Type. We showed in Lemma \ref{lem:eigenspace-decomp-array} that the portion of $\mathcal{L}$ corresponding to the bottom border of $\Delta_d$ corresponds to the $A$-eigenspace decomposition of $V$. Our next goal is to show that the portion of $\mathcal{L}$ corresponding to the left (resp. right) border of $\Delta_d$ corresponds to the $\Phi$-split (resp. $\Phi^{\Downarrow}$-split) decomposition of $V$.

For the rest of this section, fix a nonzero $\mathbf{v} \in E_0^*V$. By \cite{lp-original}, $E_i\mathbf{v} \neq 0$ for $0 \leq i \leq d$. Define $\mathbf{v}_i = E_i\mathbf{v}$ for $0 \leq i \leq d$. Then $\mathbf{v} = \mathbf{v}_0 + \cdots+ \mathbf{v}_d$ satisfies the requirements for Definition \ref{def:leonard-array} described at the beginning of Section \ref{sec:diag-type}. Recall from Definitions \ref{def:leonard-array} and \ref{def:diagonal-array} the Concrete Billiard Array $\mathcal{L}$ of Polynomial Type corresponding to $A$ and $\mathbf{v}$.

Consider again the picture from Remark $\ref{rmk:picture}$. The next two results are about the left and right borders of $\Delta_d$.
\begin{lemma} \label{lem:split-decomp}
    Let $\{U_i\}_{i=0}^d$ denote the $\Phi$-split decomposition of $V$. Consider the location $\lambda_i =(i,d-i,0)$ in $\Delta_d$. Then, $\mathcal{L}_{\lambda_i} \in U_i$.
\end{lemma}

\begin{proof}
    Use \eqref{eq:uialt} with Definition \ref{def:tau-eta} and Lemma \ref{lem:eta-zero}.
\end{proof}

\begin{lemma}
    Let $\{U^{\Downarrow}_i\}_{i=0}^d$ denote the $\Phi^{\Downarrow}$-split decomposition of $V$. Consider the location $\lambda_i =(0,d-i,i)$ in $\Delta_d$. Then, $\mathcal{L}_{\lambda_i} \in U^{\Downarrow}_i$.
\end{lemma}

\begin{proof}
    Use \eqref{eq:uialt2} with Definition \ref{def:tau-eta} and Lemma \ref{lem:eta-zero}.
\end{proof}

In \cite{lp_classification} it is shown that every Leonard system belongs to one of several families. For each family, the eigenvalues can be expressed in closed form (see \cite{lp_classification}). The most general family is called $q$-Racah. In this case, the eigenvalues have the form 
\begin{equation} \label{equ:q-racah}
    \theta_i = a + bq^i + cq^{-i} \qquad (0 \leq i \leq d)
\end{equation}
where $q,a,b,c \in \F$ with $b,c$ nonzero and $q \notin \{1,-1,0\}$. We will use formula \eqref{equ:q-racah} in the next lemma.

\begin{lemma}
    Assume the Leonard system $\Phi$ from Definition \ref{eq:ourstartingpt} has $q$-Racah type, with eigenvalue sequence (\ref{equ:q-racah}). Recall the corresponding Concrete Billiard Array $\mathcal{L}$ of Polynomial Type from above Lemma \ref{lem:split-decomp}. Then the value function for $\mathcal{L}$ from Theorem \ref{thm:la-value-function} sends 
    \begin{align*}
        (r,s,t) \mapsto q\frac{bq^{d+t-r-1} - c}{bq^{d+t-r+1} - c}
    \end{align*}
    for all $(r,s,t) \in \Delta_{d-2}$.
\end{lemma}

\begin{proof}
     Evaluate the fraction in \eqref{equ:white-value-la} using \eqref{equ:q-racah} and simplify the result.
\end{proof}

\section*{Acknowledgments}

This paper was written for an undergraduate thesis at the University of Wisconsin - Madison. The author would like to thank his advisor, Paul Terwilliger, for his advice and support.

\section*{Declaration of Competing Interest}
The author declares there are no competing interests.

\section*{Data Availability}
No data was used for the research described in this article.

\bibliographystyle{plainurl}
\bibliography{refs}

\end{document}